\documentclass[preprint,12pt]{elsarticle}

\journal{Nonlinear Analysis: Hybrid Systems}

\newif\ifreview
\newif\ifreview
\newif\ifnotreview
\reviewtrue 
\ifreview
  \notreviewfalse
\else
  \notreviewtrue
\fi

\usepackage{caption}
\usepackage{subcaption}
\usepackage{multirow}

\usepackage{algorithm}
\usepackage{algorithmic}

\usepackage{comment}

\usepackage[disable]{todonotes}

\usepackage{custom}
\renewcommand{\nodes}{modes}
\renewcommand{\node}{mode}
\renewcommand{\Csetvar}{\mathcal{S}}

\usepackage{cleveref}
\crefalias{mytheo}{theorem}
\crefalias{myprop}{proposition}
\crefalias{mydef}{definition}
\crefalias{myrem}{remark}
\crefalias{myexem}{example}
\crefalias{mylem}{lemma}
\crefalias{mycoro}{corollary}

\usepackage{mythm}

\newcommand{\tr}{\to}
\newcommand{\trs}[1][\sigma]{\tr_{#1}}

\newcommand{\Ps}{\mathcal{X}}
\newcommand{\U}{\mathcal{U}}

\newcommand{\V}{\mathcal{U}}

\newcommand{\sys}{S}

\newcommand{\cs}{\beta}
\newcommand{\signal}{signal}
\newcommand{\signals}{signals}

\newcommand{\lhcs}{HCS}
\newcommand{\lhas}{HAS}

\newcommand{\pcite}[1]{\cite{#1}}
\newcommand{\setS}{\mathcal{S}}
\newcommand{\piB}[1]{\pi_{\Image(B_{#1})^\perp}}

\newcommand{\polyset}[1]{
  \begin{subfigure}[b]{0.3\textwidth}
    \centering
      \includegraphics[trim=3.7cm 0cm 3cm 0cm, clip, width=\textwidth]{primal_#1.png}
  \end{subfigure}%
  ~~~~~~~~~~~~~~%
  \begin{subfigure}[b]{0.3\textwidth}
    \centering
      \includegraphics[trim=3.7cm 0cm 3cm 0cm, clip, width=\textwidth]{polar_#1.png}
  \end{subfigure}
}

\usepackage{siunitx}

\usepackage{tikz}
\usetikzlibrary{arrows,matrix,decorations.pathreplacing,positioning,chains,fit,shapes,calc} 
\usetikzlibrary{calc,patterns,decorations.pathmorphing,decorations.markings}

\usepackage{graphicx}

\begin{document}
\begin{frontmatter}
\title{Geometric control of hybrid systems}

\tnotetext[adhs]{This paper extends our work on continuous-time controlled invariant sets presented at ADHS 2021~\cite{legat2021geometric} to hybrid systems. Corresponding author B.~Legat.}

\author[First]{Beno\^{i}t Legat}\ead{blegat@mit.edu}
\author[Second]{Rapha\"{e}l M. Jungers}\ead{raphael.jungers@uclouvain.be}

\address[First]{LIDS,
  MIT,
  77 Massachusetts Avenue,
  Cambridge, MA 02139-4307, USA}

\address[Second]{ICTEAM,
  UCLouvain,
  4 Av. G. Lema\^itre,
  1348 Louvain-la-Neuve, Belgium}


\begin{abstract}
  \,\,\,\,In this paper, we present a geometric approach for computing
  controlled invariant sets for hybrid control systems.
  While the problem is well studied in the ellipsoidal case,
  this family is quite conservative for constrained or switched linear systems.
  We reformulate the invariance of a set as an inequality for its support function
  that is valid for any convex set.
  This produces novel algebraic conditions for the invariance of
  sets with polynomial or piecewise quadratic support function.
\end{abstract}

\end{frontmatter}

\section{Introduction}

Computing controlled invariant sets is paramount in many applications~\pcite{blanchini2015set}.
Indeed, the existence of a controlled invariant set is equivalent to the stabilizability\footnote{In the sense that the state variables can be controlled to remain bounded.} of a control system~\pcite{sontag1983lyapunov} and
a (possibly nonlinear) stabilizable state feedback can be deduced from the controlled invariant set~\pcite{barmish1985necessary}.

The stabilizability of a linear time-invariant (LTI) control system is equivalent
to the stability of its uncontrollable subspace (which is readily accessible
in its Controllability Form) \cite[Section~2.4]{wonham1974linear}.
Indeed, the eigenvalues of its controllable subspace can be fixed to any value
by a proper choice of linear state feedback.
The resulting controlled system is stable hence an invariant ellipsoid can be
determined by solving a system of linear equations~\pcite{Liapounoff1907}.
This set is also controlled invariant for the control system.
When a control system admits an ellipsoidal controlled invariant set,
it is said to be \emph{quadratically stabilizable}.
When there exists a linear state feedback such that the resulting autonomous
system admits an ellipsoidal invariant set, it is said to be
\emph{quadratically stabilizable via linear control}.

While the stabilizability of LTI control systems is equivalent to their
quadratic stabilizability via linear control,
it is no longer the case for \emph{uncertain} or \emph{switched} systems \pcite{petersen1985quadratic}.
Furthermore, it is often desirable for \emph{constrained} systems to find a controlled invariant set of maximal
volume (or which is maximal in some direction~\pcite{ahmadi2018robust}).
For such problems, the method detailed above is not suitable as it does not
take any volume consideration but more importantly, the maximal volume
invariant set may not be an ellipsoid and may not be rendered stable via
a linear control.
For this reason, a Linear Matrix Inequality (LMI) 
was devised to encapsulate the
controlled invariance of an ellipsoid via linear control \cite[Section~7.2.2]{boyd1994linear}
and the conservatism of the choice of linear control was analysed \pcite{sontag1983lyapunov}.
As the linearity of the control was found to be conservative for uncertain systems \pcite{petersen1985quadratic},
the LMI \eqref{eq:invLMI} (or \eqref{eq:invLMI_sigma} for discrete-time) was found to encapsulate controlled invariance of an ellipsoid
via \emph{any} state-feedback \cite{barmish1985necessary}.

Recent advances in control was enabled thanks to the introduction of new families of sets such as polynomial zonotopes~\cite{kochdumper2020constrained,kochdumper2021sparse}.
However, while the LMIs mentioned above have had a tremendous impact on control, the approach is limited to ellipsoids due to its algebraic nature.
An attempt to generalize it to polynomials can be found in \cite{Prajna2004}
but as detailed in \cite[Section~2]{legat2021geometric}, it is quite conservative.
The approach studied in \cite{korda2014convex} is complementary
to our method as \cite{korda2014convex} computes
outer bounds of the maximal controlled invariant sets while we compute
actual controlled invariant sets (hence inner bounds to the maximal one).

In this paper, we reinterpret the controlled invariance in a geometric/behavioural framework,
based on convex analysis,
which allows us to formulate a general condition for the controlled invariance of arbitrary
convex sets via any state-feedback in \theoref{invconvex}.
While this condition reduces to \eqref{eq:invLMI_sigma} and \eqref{eq:invLMI} for the special case of ellipsoids,
it provides a new method for computing convex controlled invariant sets with polynomial and piecewise quadratic support functions.

This paper genralizes \cite{legat2020sum}, \cite{legat2020piecewise} and \cite{legat2021geometric} into a framework for computing convex controlled invariant sets for linear hybrid control systems.
In \cite{legat2020sum}, the authors treat the particular case where the continuous dynamic at each \emph{mode} (see \cref{def:lhcs}) is trivial, i.e., $\dot{x} = 0$.
In \cite{legat2020piecewise}, the authors extends \cite{legat2020sum} to piecewise semi-ellipsoids.
In \cite{legat2021geometric}, the authors handle the particular case where there is only one mode and no \emph{transitions} (see \cref{def:lhcs}).
While \cite{legat2020sum, legat2020piecewise} covers discrete-time systems and \cite{legat2021geometric} covers continuous-time systems,
we show in this paper that the two methods can be combined to compute
controlled invariant sets for hybrid systems, exhibiting both discrete-time and continuous-time dynamics.
Using the \emph{set programming} framework~\cite{legat2020set},
this compatibility can be understood as a consequence
of the fact that the controlled invariance conditions
require the sets to be represented with their \emph{support functions} (see \cref{def:supfun})
both in discrete-time and continuous-time.

In \secref{alg}, we show how to reduce
the computation of controlled invariant sets for \emph{hybrid control systems} to the computation of
\emph{weakly invariant} sets for \emph{hybrid algebraic systems}.
In \secref{geometric}, we develop a generic condition of control invariance
for hybrid systems using our geometric approach.
We particularize it for ellipsoids (resp. sets with polynomial and piecewise quadratic support functions)
in \secref{ell} (resp. \secref{poly} and \secref{pell}).
We illustrate these new results with a numerical example in \secref{exem}.

\paragraph*{Reproducibility}
The code used to obtain the results is
published on codeocean \pcite{legat2022geometricCO}.
The set programs are reformulated
by
SetProg \pcite{legat2020set}
into
a Sum-of-Squares program
which is
reformulated into a semidefinite program
by SumOfSquares \pcite{weisser2019polynomial}
which is
solved by Mosek v8 \pcite{mosek2017mosek81034} through MathOptInterface \pcite{legat2020mathoptinterface}.


\section{Controlled invariant set}
\label{sec:alg}
In this section we define hybrid control and algebraic systems as well as the notion of invariance that will be studied in this paper.
We then show how the invariance relations between the two different classes of systems.

\begin{mydef}
  \label{def:lhcs}
  A \emph{Linear Hybrid Control System (\lhcs{})} is a system
  $\sys = (
    T,
    (A_q, B_q)_{q \in \Nodes},
    (A_\sigma, B_\sigma)_{\sigma \in \Sigma},
    (\Ps_q, \U_q)_{q \in \Nodes},
    (\V_\sigma)_{\sigma \in \Sigma}
  )$
  where $T = (\Nodes, \Sigma, \tr)$,
  $\Nodes$ is a finite set of \nodes{},
  $\Sigma$ is a finite set of \signals{} and
  $\tr \subseteq \Nodes \times \Sigma \times \Nodes$ is a set of transitions.
  We denote $(q, \sigma, q') \in \tr$ by $q \trs[\sigma] q'$.

  Given a \node{} $q \in \Nodes{}$, we denote
  the state dimension as $n_{q,x}$ and the input dimension as $n_{q,u}$
  Given a \signal{} $\sigma$, we denote the input dimension as $n_{\sigma,u}$.
  The set $\Ps_q \subseteq \mathbb{R}^{n_{q,x}}$ is the \emph{safe set}
  corresponding to \node{} $q$ and
  the sets $\U_q \subseteq \mathbb{R}^{n_{q,u}}, \V_\sigma \subseteq \mathbb{R}^{n_{\sigma,u}}$ are the sets of allowed inputs.
  For any \node{} $q$, we have
  $A_q \in \mathbb{R}^{n_{q,x} \times n_{q,x}}$,
  $B_q \in \mathbb{R}^{n_{q,x} \times n_{q,u}}$.
  For any transition $q \trs[\sigma] q'$, we have
  $A_\sigma \in \mathbb{R}^{n_{q',x} \times n_{q,x}}$,
  $B_\sigma \in \mathbb{R}^{n_{q',x} \times n_{\sigma,u}}$.

  A trajectory of $S$ is an increasing sequence of times $t_0 < t_1 < t_2 < \cdots < t_N < t_{N+1}$,
  transitions $q_0 \trs[\sigma_1] q_1 \trs[\sigma_2] \cdots \trs[\sigma_N] q_N$,
  inputs $\bar{u}_k \in \U_{\sigma_k}$ for $k \in [N]$,
  and trajectories $x_k: [t_k, t_{k+1}] \to \Ps_{q_k} \in \mathcal{C}^1$ and $u_k: [t_k, t_{k+1}] \to \U_{q_k}$ for $k = 0, 1, \ldots, N$
  satisfying:
  \begin{align*}
    \forall k \in [N], \quad x_k(t_k) & = A_{\sigma_k} x_{k-1}(t_k) + B_{\sigma_k} \bar{u}_k\\
    \forall k \in \{0, 1, \ldots, N\}, \forall t \in [t_k, t_{k+1}], \quad \dot{x}_k(t) & = A_{q_k} x_k(t) + B_{q_k} u_k(t).
  \end{align*}
\end{mydef}

The hybrid system defined in \cref{def:lhcs} may be interpreted as
a \emph{hybrid automaton}~\cite{alur1995algorithmic}
where the guard of each transition $q \trs q'$ is $\Ps_q$ or $\R^{n_{q,x}}$.
In this context, the discrete-time dynamical system $x^+ = A_\sigma x + B_\sigma u$ is commonly referred to as the \emph{reset map}.
We allow the state space of different \nodes{} to differ as our method naturally
extends to different state spaces but the reader may consider them to have identical dimension for simplicity.

\begin{mydef}[Controlled invariant sets for a \lhcs{}]
  \label{def:cis}
  Consider a \lhcs{} $\sys$ as defined in \cref{def:lhcs}.
  We say that closed sets $\Csetvar_q \subseteq \Ps_q$ for $q \in \Nodes$ are
  \emph{controlled invariant} for $\sys$ if
  \begin{align}
    \label{eq:cis_sigma}
    \forall q \trs q', \forall x \in \Csetvar_q,
    \exists u \in \U_\sigma \text{ such that}
    \quad
    A_\sigma x + B_\sigma u & \in \Csetvar_{q'}\\
    \label{eq:cis_q}
    \forall q \in \Nodes, \forall x \in \Csetvar_q, \exists u \in \U_q \text{ such that}
    \quad
    A_q x + B_q u & \in T_{\Csetvar_{q}}(x).
  \end{align}
\end{mydef}

\Cref{eq:cis_q} is commonly referred to as the \emph{Nagumo condition}; see \cite[Theorem~4.7]{blanchini2015set}.
In view of \cref{def:cis}, the transitions are considered
\emph{autonomous} and not \emph{controlled}; see details in \cite[Section~1.1.3]{liberzon2012switching}.

\subsection{Linear Hybrid Algebraic System}
\label{sec:lhas}

In this section, we show the equivalence of the notion of invariance with
another class of systems that directly models
the geometric behaviours of the trajectories of
a \lhcs{} with unconstrained input.
The reduction of the computation of controlled invariant
sets of \lhcs{} with constrained input to \lhcs{} of unconstrained input is detailed in \cite[Section~2.2]{legat2020sum}.

\begin{mydef}
  \label{def:lhas}
  A \emph{Linear Hybrid Algebraic System (\lhas{})} is a system
  $\sys = (
    T,
    (C_q, E_q)_{q \in \Nodes},
    (C_\sigma, E_\sigma)_{\sigma \in \Sigma},
    (\Ps_q)_{q \in \Nodes}
  )$
  where $T = (\Nodes, \Sigma, \tr)$,
  $\Nodes$ is a finite set of \nodes{},
  $\Sigma$ is a finite set of \signals{} and
  $\tr \subseteq \Nodes \times \Sigma \times \Nodes$ is a set of transitions.

  Given a \node{} $q \in \Nodes{}$, we denote
  the state dimension as $n_{q,x}$.
  The set $\Ps_q \subseteq \mathbb{R}^{n_{q,x}}$ is the \emph{safe set}
  corresponding to \node{} $q$.
  For any \node{} $q$,
  there exists a $n_{q,p}$ such that,
  $C_q \in \mathbb{R}^{n_{q,p} \times n_{q,x}}$,
  $E_q \in \mathbb{R}^{n_{q,p} \times n_{q,x}}$.
  For any transition $q \trs[\sigma] q'$,
  there exists a $n_{\sigma,p}$ such that,
  $C_\sigma \in \mathbb{R}^{n_{\sigma,p} \times n_{q,x}}$,
  $E_\sigma \in \mathbb{R}^{n_{\sigma,p} \times n_{q',x}}$.

  A trajectory of $S$ is an increasing sequence of
  times $t_0 < t_1 < t_2 < \cdots t_N$,
  transitions $q_0 \trs[\sigma_1] q_1 \trs[\sigma_2] \cdots \trs[\sigma_N] q_N$,
  and trajectories $x_k: [t_{k-1}, t_k] \to \Ps_{q_k} \in \mathcal{C}^1$ for $k = 0, 1, \ldots, N$
  satisfying:
  \begin{align*}
    \forall k \in [N], \quad E_{\sigma_k} x_k(t_k) & = C_{\sigma_k} x_{k-1}(t_k)\\
    \forall k \in \{0, 1, \ldots, N\}, \forall t \in [t_{k-1}, t_k], \quad E_{q_k} \dot{x}_k(t) & = C_{q_k} x_k(t).
  \end{align*}
\end{mydef}

\begin{mydef}[Weakly invariant sets for a \lhas{}]
  \label{def:is}
  Consider a \lhas{} $\sys$ as defined in \cref{def:lhas}.
  We say that closed sets $\Csetvar_q \subseteq \Ps_q$ for $q \in \Nodes$ are
  \emph{weakly invariant} for $\sys$ if
  \begin{align}
    \label{eq:is_sigma}
    \forall q \trs q', \forall x \in \Csetvar_q,
    \quad
    C_\sigma x & \in E_\sigma \Csetvar_{q'}\\
    \label{eq:is_q}
    \forall q \in \Nodes, \forall x \in \Csetvar_q,
    \quad
    C_q x & \in E_q T_{\Csetvar_{q}}(x).
  \end{align}
\end{mydef}



We now show that the computation of controlled invariant sets for a \lhcs{} can be reduced to
the computation of weakly invariant sets for a \lhas{}.

\begin{mylem}[{\cite[Proposition~4]{legat2021geometric}}]
  \label{lem:projB}
  Given a subset $\mathcal{S} \subseteq \R^n$ and
  matrices $A \in \R^{r \times n}, B \in \R^{r \times m}$, the following holds:
  \[ A \mathcal{S} + B \R^m = \piB{}^{-1} \piB{} A \mathcal{S} \]
  where $\piB{}$ is any orthogonal projection matrix onto the orthogonal subspace of $\Image(B)$
  and $\piB{}^{-1}$ is the preimage defined in \cref{eq:preimage}.
  \begin{proof}
    Given $x \in \mathcal{S}$ and $y \in \R^r$, we have $y \in A \{x\} + B \R^m$ if and only if $y - Ax \in \Image(B)$.
    As $\piB{}$ is orthogonal, its kernel is $\Image(B)$.
    Therefore $y - Ax \in \Image(B)$ is equivalent to $\piB y = \piB Ax$.
  \end{proof}
\end{mylem}

The following proposition generalizes both \cite[Proposition~2]{legat2020sum} and \cite[Proposition~5]{legat2021geometric}.

\begin{myprop}
  \label{prop:proju}
  The sets $\Csetvar = (\Csetvar_q)_{q \in \Nodes}$ are \emph{controlled invariant} for the \lhcs{}
  $\sys = (
    T,
    (A_q, B_q)_{q \in \Nodes},
    (A_\sigma, B_\sigma)_{\sigma \in \Sigma},
    (\Ps_q, \R^{n_{q,u}})_{q \in \Nodes},
    (\R^{n_{\sigma,u}})_{\sigma \in \Sigma}
  )$
  if and only if
  they are weakly invariant sets for the \lhas{}
  $$\sys' = (
    T,
    (\piB{q} A_q, \piB{q})_{q \in \Nodes},
    (\piB{\sigma} A_\sigma, \piB{\sigma})_{\sigma \in \Sigma},
    (\Ps_q)_{q \in \Nodes}
  ).$$

  \begin{proof}
    By \cref{lem:projB}, \eqref{eq:cis_sigma} is equivalent
    \[ \piB{\sigma} A_\sigma x \in \piB{\sigma} \Csetvar_{q'}. \]
    which is \eqref{eq:is_sigma} for $\sys'$.

    Similarly, by \cref{lem:projB}, \eqref{eq:cis_q} is equivalent
    \[ \piB{q} A_q x \in \piB{q} T_{\Csetvar_{q}}(x). \]
    which is \eqref{eq:is_q} for $\sys'$.
  \end{proof}
\end{myprop}

\section{Computing controlled invariant sets}
\label{sec:geometric}

In this section we derive a characterization of
the weak invariance of closed convex sets
under the form of inequalities for their support functions.
This section uses notions of convex geometry
that are recalled in \ref{sec:convex}.
The following theorem generalizes both \cite[(27)]{legat2020piecewise} and \cite[Theorem~7]{legat2021geometric}.

\begin{mytheo} 
  \label{theo:invconvex}
  Consider a \lhas{} $\sys$ as defined in \cref{def:lhas}.
  Closed sets $\Csetvar_q \subseteq \Ps_q$ for $q \in \Nodes$ are
  weakly invariant for $\sys$
  if and only if
  \begin{align}
    \label{eq:inv_sigma}
    \forall q \trs q', \forall y \in \R^{n_{\sigma,p}}, \quad \supfun{C_\sigma^\Tr y}{\Csetvar_q} \le \supfun{E_\sigma^\Tr y}{\Csetvar_{q'}}\\
    \label{eq:invface}
    \forall q \in \Nodes, \forall z \in \mathbb{R}^{n_{q,p}}, \forall x \in F_\Csetvar(E_q^\Tr z), \quad \la z, C_q x \ra \le 0
  \end{align}
  where $F_\Csetvar$ denotes the \emph{exposed face} defined in \cref{def:exposed}.
  \begin{proof}
    We start by proving the equivalence between \eqref{eq:is_sigma} and \eqref{eq:inv_sigma}.
    By \cref{prop:supportincl}, \cref{eq:is_sigma} is equivalent to
    \[ \forall q \trs q', \forall y \in \R^{n_{\sigma,p}}, \quad \supfun{y}{C_\sigma \Csetvar_q} \le \supfun{y}{E_\sigma \Csetvar_{q'}} \]
    which is equivalent to \cref{eq:inv_sigma} by \cref{prop:hA}.

    We now prove the equivalence between \eqref{eq:is_q} and \eqref{eq:invface}.
    Given any \node{} $q$, as $\Csetvar_q$ is convex, $T_{\Csetvar_q}(x)$ is a convex cone. By definition
    of the polar of a cone, $x \in E_q T_{\Csetvar_q}(x)$ if and only if
    $\langle y, x \rangle \le 0$ for all $y \in \polar{[E_q T_{\Csetvar_q}(x)]}$.
    By \propref{podu}, $\polar{[E_q T_{\Csetvar_q}(x)]} = E_q^{-\Tr}N_{\Csetvar_q}(x)$.
    Therefore, the set $\Csetvar_q$ is weakly invariant if and only if
    \begin{equation*}
      \forall x \in \partial \Csetvar_q, \forall z \in E_q^{-\Tr} N_{\Csetvar_q}(x), \langle z, C_q x \rangle \le 0.
    \end{equation*}
    By \propref{normalsupfun}, we have
    \begin{multline*}
      \{\, (x, z) \in \partial \Csetvar_q \times \mathbb{R}^r \mid E^\Tr z \in N_{\Csetvar_q}(x) \,\} = \\
      \{\, (x, z) \in \partial \Csetvar_q \times \mathbb{R}^r \mid x \in F_{\Csetvar_q}(E_q^\Tr z) \,\}.
    \end{multline*}
  \end{proof}
\end{mytheo}

As we show in the remainder of this section, \cref{theo:invconvex} allows to reformulate
the invariance as an inequality in terms of the support functions of the sets $\Csetvar_q$.
This is already the case of \cref{eq:inv_sigma} so it remains to reformulate \cref{eq:invface}. 
As shown in the following theorem, this is possible in case the support function is differentiable.
We generalize this result with a relaxed notion of differentiability in \theoref{pdiff}.
The following theorem generalizes both \cite[(27)]{legat2020piecewise} and \cite[Theorem~8]{legat2021geometric}.

\begin{mytheo}
  \label{theo:invexposed}
  Consider a \lhas{} $\sys$ as defined in \cref{def:lhas} and
  nonempty closed convex sets $\Csetvar_q \subseteq \Ps_q$ for $q \in \Nodes$
  such that
  $\supfun{\cdot}{\Csetvar_q}$ is differentiable for all $q \in \Nodes$.
  Then the sets are weakly invariant for $\sys$
  if and only if
  \begin{align}
    \notag
    \forall q \trs q', \forall y \in \R^{n_{\sigma,p}}, \supfun{C_\sigma^\Tr y}{\Csetvar_q} & \le \supfun{E_\sigma^\Tr y}{\Csetvar_{q'}}\\
    \label{eq:invgrad}
    \forall q \in \Nodes, \forall z \in \mathbb{R}^{n_{q,p}}, \la z, C_q \grad \supfun{E_q^\Tr z}{\Csetvar_q} \ra & \le 0.
  \end{align}
  \begin{proof}
    By \propref{exposed}, $F_{\Csetvar_q}(E_q^\Tr z) = \{\grad \supfun{E_q^\Tr z}{\Csetvar_q}\}$
    hence \eqref{eq:invface} is equivalent to \eqref{eq:invgrad}.
  \end{proof}
\end{mytheo}

As \cref{theo:invexposed} formulates the invariance in terms of the support function of $\Csetvar_q$,
it allows to combine the invariance constraint with other set constraints that can be formulated in terms of support functions.
Moreover, for an appropriate family of sets, also called \emph{template}, the set program
can be automatically rewritten into a convex program combining all constraints using \emph{set programming} \cite{legat2020set}.
For this reason, we only focus on the invariance constraint and do not detail how to
formulate the complete convex programs with the objective and all the constraints needed to obtain the results of \secref{exem} as these problems are decoupled.

\subsection{Ellipsoidal controlled invariant set}
\label{sec:ell}



In this section, we particularize \theoref{invexposed} to the case of ellipsoids.
Since the support function of an ellipsoid $\mathcal{E}_Q$ is $\supfun{y}{\mathcal{E}_Q} = \sqrt{y^\Tr Q^{-1} y}$,
we have the following corollary of \theoref{invexposed}
that generalizes both \cite[Theorem~2]{legat2020sum} and \cite[Corollary~9]{legat2021geometric}.
\begin{mycoro}
  Consider a \lhas{} $\sys$ as defined in \cref{def:lhas}
  and positive semidefinite matrices $Q_q$ such that
  the ellipsoid $\mathcal{E}_Q \subseteq \Ps_q$ for $q \in \Nodes$.
  Then the sets are weakly invariant for $\sys$
  if and only if
  \begin{align}
    \label{eq:invLMI_sigma}
    \forall q \trs q', C_\sigma Q_q^{-1} C_\sigma^\Tr \preceq E_\sigma Q_{q'}^{-1} E_\sigma^\Tr\\
    \label{eq:invLMI}
    \forall q \in \Nodes, C_qQ_q^{-1}E_q^\Tr + E_qQ_q^{-1}C_q^\Tr \preceq 0.
  \end{align}
\end{mycoro}

Observe that for the trivial case $\Image(B_q) = \mathbb{R}^{n_{q,x}}$ for some node $q$,
\cref{prop:proju} produces a \lhas{} with $n_{q,p} = 0$ hence the
LMI~\eqref{eq:invLMI} will be trivially satisfied for any $Q_q^{-1}$, which is expected.
The same applies for \eqref{eq:invLMI_sigma} in case $\Image(B_\sigma) = \mathbb{R}^{n_{\sigma,x}}$ for some $\sigma$.

\subsection{Polynomial controlled invariant set}
\label{sec:poly}
In this section, we derive the algebraic condition for
the controlled invariance of a set with polynomial support function.
This template is referred to as \emph{polyset}; see \cite[Section~1.5.3]{legat2020set}.
The following corollary generalizes both \cite[Theorem~5]{legat2020sum} and \cite[Corollary~10]{legat2021geometric}.
\begin{mycoro}
  \label{coro:poly}
  Consider a \lhas{} $\sys$ as defined in \cref{def:lhas},
  convex homogeneous\footnote{A polynomial is \emph{homogeneous} if all its monomials have the same total degree} nonnegative polynomials $(p_q(x))_{q \in \Nodes}$ of degree $2d$
  and
  the sets $\Csetvar_q$
  defined by the support function
  $\supfun{y}{\Csetvar_q} = p_q(y)^{\frac{1}{2d}}$ for $q \in \Nodes$.
  Suppose that $\Csetvar_q \subseteq \Ps_q$ for all $q \in \Nodes$.
  Then the sets are weakly invariant for $\sys$
  if and only if
  \begin{align}
    \label{eq:invSOS_sigma}
    \forall q \trs q', \forall y \in \R^{n_{\sigma,p}}, p_q(C_\sigma^\Tr y) \le p_{q'}(E_\sigma^\Tr y)\\
    \label{eq:invSOS}
    \forall q \in \Nodes, \forall z \in \R^{n_{q,p}}, z^\Tr C_q \grad p_q(E_q^\Tr z) \le 0.
  \end{align}
  \begin{proof}
    We have
    \[
      \grad \supfun{y}{\Csetvar_q} = \frac{1}{p_q(y)^{1-\frac{1}{2d}}} \grad p_q(y).
    \]
    If $p(y)$ is identically zero, this is trivially satisfied.
    Otherwise, $p_q(y)^{1-\frac{1}{2d}}$ is nonnegative and is zero
    in an algebraic variety of dimension $n - 1$ at most.
    Therefore,
    \eqref{eq:invgrad} is equivalent to \eqref{eq:invSOS}.
  \end{proof}
\end{mycoro}

The conditions \eqref{eq:invSOS_sigma} and \eqref{eq:invSOS}
require the nonnegativity of a multivariate polynomial.
While verifying the nonnegativity of a polynomial is co-NP-hard,
a sufficient condition can be obtained via the standard Sum-of-Squares programming framework; see~\cite{blekherman2012semidefinite}.
Moreover, the theorem requires the convexity of the polynomials $p_q$.
It is shown in~\cite{ahmadi2013np} that the convexity or quasi-convexity of a multivariate polynomial of degree at least four is NP-hard to decide.
However, the convexity constraint can be replaced by the tractable SOS-convexity
constraint which is a sufficient condition for convexity~\cite{ahmadi2013np}.

\subsection{Piecewise semi-ellipsoidal controlled invariant set}
\label{sec:pell}

In \cite{johansson1998computation}, the authors study the computation of piecewise quadratic Lyapunov functions for continuous-time autonomous piecewise affine systems.
In \cite{legat2020piecewise}, the authors present a convex programming approach
to compute \emph{piecewise semi-ellipsoidal} controlled invariant sets for discrete-time control systems.
A similar approach is developed in \cite{legat2021geometric} for continuous-time control system.
In this section, we combine the two approaches into a condition for hybrid systems using \theoref{invconvex}.
We recall \cite[Definition~2]{legat2020piecewise} below.
\begin{mydef}
  A \emph{polyhedral conic partition} of $\R^n$ is a set of $m$ polyhedral cones $(\mathcal{P}_i)_{i=1}^m$
  with nonempty interior
  such that for all $i \neq j$, $\dim(\mathcal{P}_i \cap \mathcal{P}_j) < n$
  and $\cup_{i=1}^m \mathcal{P}_i = \R^n$.
\end{mydef}
A polyhedral conic partition defines the full-dimensional faces of a \emph{complete fan}, as defined in \cite[Section~7]{ziegler1995lectures}.

A piecewise semi-ellipsoid has a support function of the form
\begin{equation}
  \supfun{y}{\Csetvar} = \sqrt{y^\Tr Q_i y}, \qquad y \in \mathcal{P}_i, \qquad i = 1, \ldots, m
\end{equation}
where $(\mathcal{P}_i)_{i=1}^m$ is a polyhedral conic partition.
The support function additionally has to satisfy \cite[(2) and (3)]{legat2020piecewise}
to ensure its continuity and convexity.
The following theorem generalizes both \cite[(27)]{legat2020piecewise} and \cite[Theorem~12]{legat2021geometric}.

\begin{mytheo}
  \label{theo:pdiff}
  Consider a \lhas{} $\sys$ as defined in \cref{def:lhas},
  polyhedral conic partitions $(\mathcal{P}_{q,i})_{i=1}^{m_q}$ and
  nonempty closed convex sets $(\Csetvar_q)_{q \in \Nodes}$ defined by the support function
  \[ \supfun{y}{\Csetvar_q} = f_{q,i}(y) \qquad y \in \mathcal{P}_{q,i} \qquad i = 1, \ldots, m_q. \]
  Suppose that $\Csetvar_q \subseteq \Ps_q$ for all $q \in \Nodes$.
  The sets $\Csetvar_q$ are weakly invariant for $\sys$
  if and only if
  \begin{align}
    \notag
    \forall q \trs q', \forall i \in [m_q], \forall j \in [m_{q'}],\qquad\qquad\qquad\qquad\\
    \label{eq:pdiff_sigma}
    \qquad \forall y \in C_\sigma^{-\Tr} \mathcal{P}_{q,i} \cap E_\sigma^{-\Tr} \mathcal{P}_{q',j}, \quad f_{q,i}(C_\sigma^\Tr y) & \le f_{q',j}(E_\sigma^\Tr y)\\
    \label{eq:invpgrad}
    \forall q \in \Nodes, \forall i \in [m_q], \forall z \in E_q^{-\Tr} \mathcal{P}_{q,i}, \quad \la z, C_q \grad f_{q,i}(E_q^\Tr z) \ra & \le 0.
  \end{align}
  \begin{proof}
    If $y \in C_\sigma^{-\Tr} \mathcal{P}_{q,i} \cap E_\sigma^{-\Tr} \mathcal{P}_{q',j}$, then $\supfun{C_\sigma^\Tr y}{\Csetvar_q} = f_{q,i}(C_\sigma^\Tr y)$ and $\supfun{E_\sigma^\Tr y}{\Csetvar_{q'}} = f_{q',j}(E_\sigma^\Tr y)$
    hence \eqref{eq:inv_sigma} is reformulated as \eqref{eq:pdiff_sigma}.

    We now prove the equivalence between \eqref{eq:invface} and \eqref{eq:invpgrad}.
    Consider a \node{} $q \in \Nodes$.
    Given $z \in \mathbb{R}^{n_{q,p}}$ such that $E_q^\Tr z$ is in the intersection of the boundary of $\Csetvar_q$ and
    the interior of $\mathcal{P}_{q,i}$,
    the support function is differentiable at $E_q^\Tr z$
    hence, by \propref{exposed},
    $F_\Csetvar(E_q^\Tr z) = \{\grad f_{q,i}(E_q^\Tr z)\}$.
    The condition~\eqref{eq:invface} is therefore reformulated as \eqref{eq:invpgrad}.

    Given a subset $I$ of $\{1, \ldots, m\}$ and
    $z \in \mathbb{R}^{n_{q,p}}$ such that $E_q^\Tr z$ is in
    the intersection of the boundary of $\Csetvar_q$ and $\cap_{i \in I} \mathcal{P}_{q,i}$,
    $F_{\Csetvar_q}(E_q^\Tr z)$ is the convex hull of
    $\grad \supfun{E_q^\Tr z}{\Csetvar_q}$ for all $i \in I$.
    For any convex combination (i.e., nonnegative numbers summing to 1) $(\lambda_i)_{i \in I}$,
    \eqref{eq:invpgrad} implies that
    \[ \la z, C_q \sum_{i \in I} \lambda_i \grad f_{q,i}(E_q^\Tr z) \ra = \sum_{i \in I} \lambda_i \la z, C_q \grad f_{q,i}(E_q^\Tr z) \ra \le 0. \]
  \end{proof}
\end{mytheo}

The following corollary generalizes both \cite[Theorem~4]{legat2020piecewise} and \cite[Corollary~13]{legat2021geometric}.
\begin{mycoro}
  Consider a \lhas{} $\sys$ as defined in \cref{def:lhas} and
  piecewise semi-ellipsoids $\Csetvar_q \subseteq \Ps_q$ for $q \in \Nodes$.
  The sets are weakly invariant for $\sys$
  if and only if
  \begin{align}
    \notag
    \forall q \trs q', \forall i \in [m_q], \forall j \in [m_{q'}],\qquad\qquad\qquad\qquad\qquad\qquad\qquad\quad\\
    \label{eq:invpell_sigma}
    \qquad \forall y \in C_\sigma^{-\Tr} \mathcal{P}_{q,i} \cap E_\sigma^{-\Tr} \mathcal{P}_{q',j}, y^\Tr E_\sigma Q_{q',j} E_\sigma^\Tr y \le y^\Tr E_\sigma Q_{q',j} E_\sigma^\Tr y\\
    \label{eq:invpell}
    \forall q \in \Nodes, \forall i \in [m_q], \forall z \in E_q^{-\Tr} \mathcal{P}_{q,i}, \quad z^\Tr C_q Q_{q,i}^{-1}E_q^\Tr z + z^\Tr E_q Q_{q,i}^{-1}C_q^\Tr z \le 0.
  \end{align}
\end{mycoro}

The conditions \eqref{eq:invpell_sigma} and \eqref{eq:invpell} amount to verifying the positive semidefiniteness of a quadratic form
when restricted to a polyhedral cone. When this cone is the positive orthant, this is called the
\emph{copositivity} which is co-NP-complete to decide \pcite{murty1987some}.
However, a sufficient LMI is given in \cite[Proposition~2]{legat2020piecewise}
and a necessary and sufficient condition is given by a hierarchy of Sum-of-Squares programs \cite[Chapter~5]{parrilo2000structured}.
We use the sufficient LMI in the numerical example of \cref{sec:exem}.

\section{Numerical example}
\label{sec:exem}

This example considers the \lhcs{} with one \node{} of continuous-time dynamics:
\begin{align*}
\dot{x}_1(t) & = x_2(t)\\
\dot{x}_2(t) & = u(t)
\end{align*}
with state constraint $x \in [-1, 1]^2$ and input constraint $u \in [-1, 1]$
and the following transition from the only \node{} to itself:
\begin{align*}
x_1^+ & = -x_1 + u/8\\
x_2^+ & = x_2 - u/8
\end{align*}
with state constraint $x \in [-1, 1]^2$ and input constraint $u \in [-1, 1]$.

The union of controlled invariant sets is controlled invariant.
Moreover, by linearity and convexity of the constraint sets, the convex hull of the unions of controlled
invariant sets is controlled invariant.
Therefore, there exists a \emph{maximal} controlled invariant set, i.e., a controlled invariant set in which all controlled invariant sets are included, for any family that is closed under union (resp. convex hull);
it is the union (resp. convex hull) of all controlled invariant sets included in $[-1, 1]^2$.

\begin{figure}[!ht]
  \centering
  \polyset{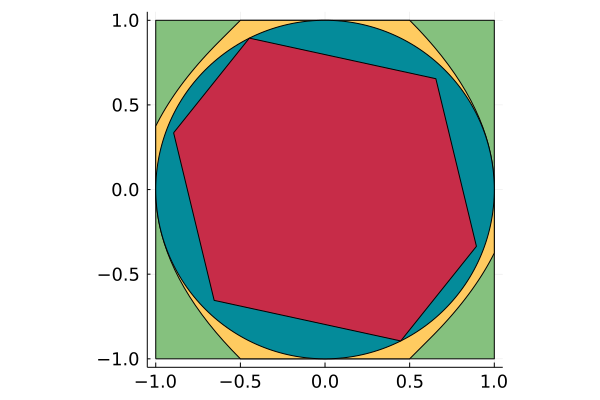}
  \polyset{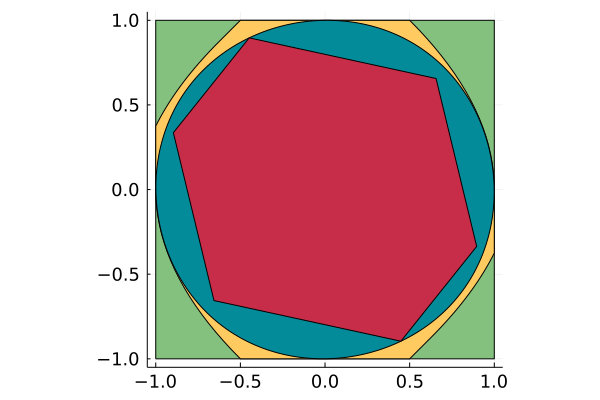}
  \polyset{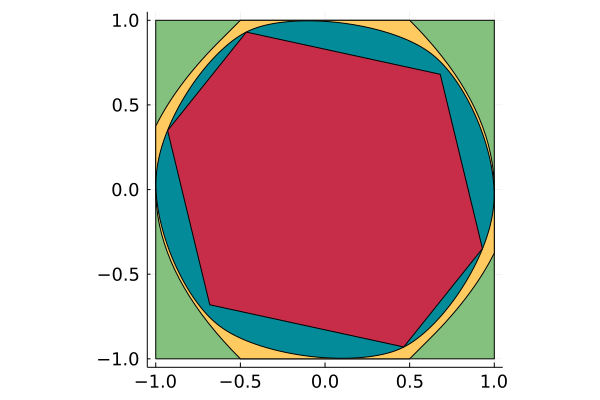}
  \polyset{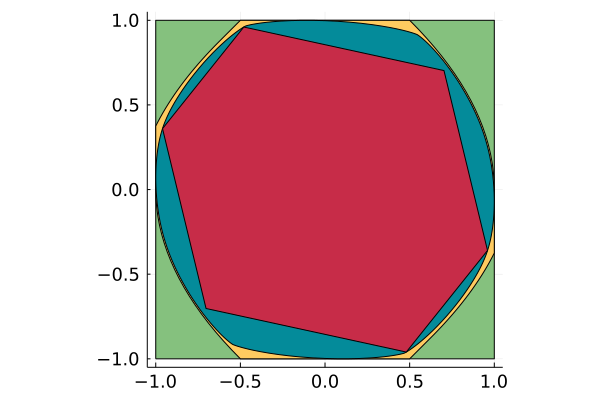}
  \caption{
    In blue are the solution for polysets of different degrees.
    The degrees from top to bottom are respectively 2, 4, 6 and 8.
    The green set is the safe set $[-1, 1]^2$, the yellow set is the maximal
    controlled invariant set and the red set is $\gamma \mathcal{D}$.
    The sets are represented in the primal space in left figures
    and in polar space in the right figures.
  }
  \label{fig:poly}
\end{figure}

\begin{figure}[!ht]
  \centering
  \polyset{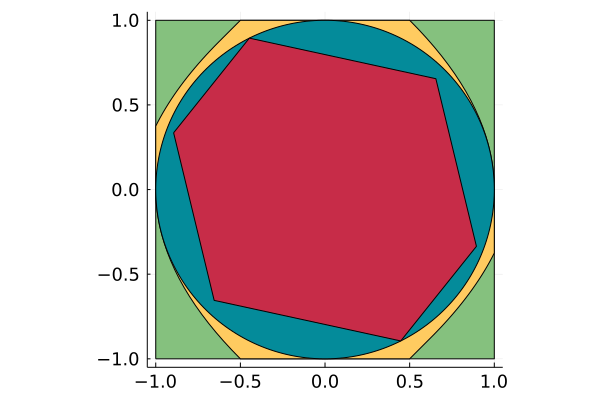}
  \polyset{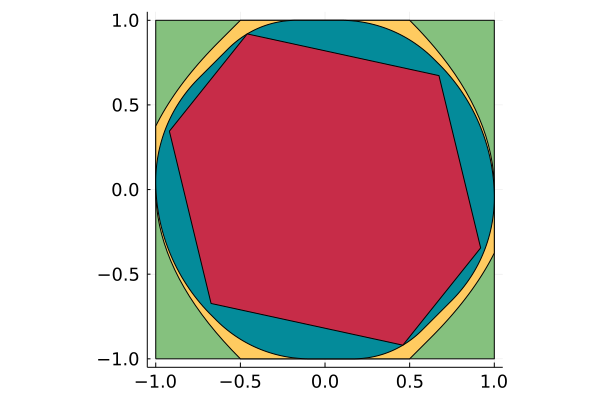}
  \polyset{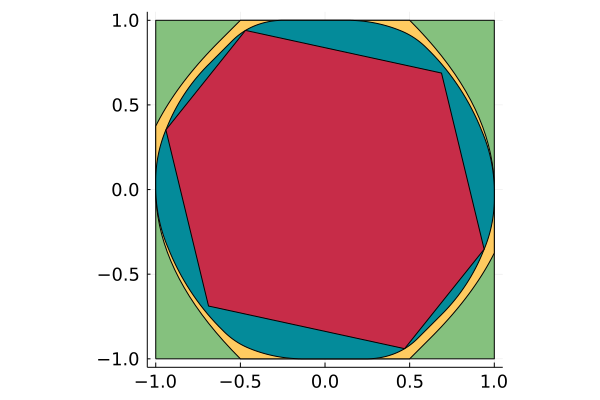}
  \caption{
    In blue are the solution for piecewise semi-ellipsoids for two
    different polyhedral conic partitions.
    The partitions from top to bottom are as described in \eqref{eq:sphere}
    with $m = (4, 3)$ (resp. $(8, 5)$, $(16, 7)$).
    The green set is the safe set $[-1, 1]^2$, the yellow set is the maximal
    controlled invariant set and the red set is $\gamma \mathcal{D}$.
    The sets are represented in the primal space in left figures
    and in polar space in the right figures.
  }
  \label{fig:pell}
\end{figure}

For this simple planar system, the maximal controlled invariant set
can be obtained by hand. 
We represent it in yellow in \figref{poly} and \figref{pell}.

As \cref{prop:proju} requires the input to be unconstrained,
it cannot be applied to this system directly.
We follow the approach detailed in \cite[Section~2.2]{legat2020sum}
to reduce the computation of controlled invariant sets for this system to a system with uncontrolled input.
In this example, it corresponds to
the projection onto the first two dimensions of controlled invariant sets for the
following lifted system:
\begin{align*}
\dot{x}_1(t) & = x_2(t)\\
\dot{x}_2(t) & = x_3(t)\\
\dot{x}_3(t) & = u(t)
\end{align*}
with state constraint $x \in [-1, 1]^3$;
with a first transition to a temporary \node{}:
\begin{align*}
x_1^+ & = x_1\\
x_2^+ & = x_2\\
x_3^+ & = u
\end{align*}
with state constraint $x \in [-1, 1]^3$ and unconstrained input;
and a second transition back to the original \node{}:
\begin{align*}
x_1^+ & = -x_1 + x_3/8\\
x_2^+ & = x_2 - x_3/8\\
x_3^+ & = u.
\end{align*}
Note that the input $u$ chosen in the first transition is the input that will be used for
the reset map and the input $u$ chosen for the second jump is the input that will be used
for the state $x_3$ of the continuous-time system.

As shown in \cref{prop:proju}, a set is controlled invariant for this system if and only if it is weakly invariant for the algebraic system
\begin{align*}
\dot{x}_1(t) & = x_2(t)\\
\dot{x}_2(t) & = x_3(t)
\end{align*}
with state constraint $x \in [-1, 1]^3$;
with a first transition to a temporary \node{}:
\begin{align*}
x_1^+ & = x_1\\
x_2^+ & = x_2
\end{align*}
with state constraint $x \in [-1, 1]^3$
and a second transition back to the original \node{}:
\begin{align*}
x_1^+ & = -x_1 + x_3/8\\
x_2^+ & = x_2 - x_3/8.
\end{align*}

We represent the state set $[-1, 1]^2$ and its polar in green in \figref{poly} and \figref{pell}.

While the \emph{maximal} invariant set is well defined,
it is not the case anymore when we restrict the set to belong to the family of ellipsoids,
polysets or piecewise semi-ellipsoids for a fixed polyhedral conic partition
as these families are not invariant under union nor convex hull.
The objective used to determine which invariant set is selected
depends on the particular application.
Let $\mathcal{D}$ be the convex hull of $\{(-1 + \sqrt3, -1 + \sqrt3), (-1/2, 1), (-1, 3/4), (1 - \sqrt3, 1 - \sqrt3), (1/2, -1), (1, -3/4)\}$.
For this example, we maximize $\gamma$ such that $\gamma\mathcal{D}$ is
included in the projection of the invariant set onto the first two dimensions.
We represent $\gamma\mathcal{D}$ in red in \figref{poly} and \figref{pell}.

For the ellipsoidal template considered in \secref{ell}, the optimal solution
is shown in \figref{poly} as ellipsoids corresponds to polysets of degree 2.
The optimal objective value is $\gamma \approx 0.894$.

For the polyset template considered in \secref{poly}, the optimal solution
are represented in \figref{poly}. 
The optimal objective value for degree 4 (resp. 6 and 8)
is $\gamma \approx 0.896$.
(resp.
$\gamma \approx 0.93$ and
$\gamma \approx 0.96$).

For the piecewise semi-ellipsoidal template, we consider
as polyhedral conic partitions
the \emph{face fan}~\cite[Example~7.2]{ziegler1995lectures}, i.e., the conic hull of each facet, of the polytope with extreme points
\begin{equation}
  \label{eq:sphere}
  (\cos(\alpha)\cos(\beta), \sin(\alpha)\cos(\beta), \sin(\beta))
\end{equation}
where $\alpha = 0, 2\pi/m_1, 4\pi/m_1, \ldots, 2(m_1-1)\pi/m_1$ and
$\beta = -\pi/2, \ldots, -2\pi/(m_2-1), -\pi/(m_2-1), 0, \pi/(m_2-1), 2\pi/(m_2-1), \ldots, \pi/2$.

The optimal objective value for $m = (4, 3)$ (resp. $(8, 5)$, $(16, 7)$) is $\gamma \approx 0.894$ (resp. $\gamma \approx 0.92$, $\gamma \approx 0.94$).
The corresponding optimal solution is shown in \figref{pell}.

\section{Conclusion}
We proved a condition for controlled invariance of convex sets for a hybrid control system
based on their support functions.
We particularized the condition for three templates:
ellipsoids, polysets and piecewise semi-ellipsoids.
In the ellipsoidal case, it combines known LMIs for discrete-time and continuous-time systems.
In the polyset case, it provides a condition significantly less conservative than~\cite{Prajna2004}.
Indeed, our condition is equivalent to invariance by \cororef{poly} and, as shown in \cite[Section~2]{legat2021geometric}, \cite{Prajna2004} is quite conservative.
In the piecewise semi-ellipsoidal case, it provides the first convex programming approach for the controlled invariance of hybrid control systems to the best of our knowledge.

As future work, we aim to apply this framework to other families such as
the \emph{piecewise polysets} defined in \cite{legat2020set}.
Moreover, instead of considering a uniform discretization of the hypersphere
as in \eqref{eq:sphere}, a more adaptive methods could be considered.
The sensitivity information provided by the dual solution of the optimization program could for instance determine
which pieces of the partition should be refined.

Finally, our definition of hybrid control system (\cref{def:lhcs}) does not support encoding a \emph{guard}
that would restrict the possible transitions depending on the current state.
Integrating this additional feature to the framework would allow the method to handle \emph{any} hybrid automaton
with linear continuous-time dynamic at each mode and linear reset maps.



\bibliographystyle{plain}
\bibliography{biblio}

\section*{Acknowledgements}

The first author is a post-doctoral fellow of the Belgian American Educational Foundation.
His work is partially supported by the National Science Foundation under Grant No. OAC-1835443.
The second author is a FNRS honorary Research Associate.
This project has received funding from the European Research Council (ERC) under the European Union's Horizon 2020 research and innovation programme under grant agreement No 864017 - L2C. RJ is also supported by the Innoviris Foundation and the FNRS (Chist-Era Druid-net)

\appendix
\section{Convex geometry}
\label{sec:convex}

\begin{mydef}[{
    \cite[p.~28]{rockafellar2015convex}}]
  \label{def:supfun}
  Consider a convex set $\Csetvar$.
  The \emph{support function} of $\Csetvar$ is defined as
\[
  \supfun{y}{\Csetvar} = \sup_{x \in \Csetvar} \la y, x \ra.
\]
\end{mydef}

\begin{mydef}[{Polar set}]
  For any convex set $\Csetvar$ the polar of $\Csetvar$,
  denoted $\polar{\Csetvar}$,
  is defined as
\[
    \polar{\Csetvar} = \setbuild{y}{\supfun{y}{\Csetvar} \le 1}.
\]
\end{mydef}

We define the \emph{tangent cone} as follows \cite[Definition~4.6]{blanchini2015set}.

\begin{mydef}[{Tangent cone}]
  Given a closed convex set $\setS$ and a distance function $\dist{\setS}{x}$,
  the \emph{tangent cone} to $\setS$ at $x$ is defined as follows:
  \[ T_{\setS}(x) = \left\{\, y \mid \lim_{\tau \to 0} \frac{\dist{\setS}{x + \tau y}}{\tau} = 0 \,\right\} \]
  where the distance is defined as
\[
    \dist{\setS}{x} = \inf_{y \in \setS} \|x - y\|
\]
  where $\|\cdot\|$ is a norm.
  The tangent cone is a convex cone and is independent of the norm used.

\end{mydef}

For a convex set $\Csetvar$, the \emph{normal cone} is the polar of the tangent cone $N_\Csetvar(x) = \polar{T_\Csetvar}(x)$.

The \emph{exposed face} (also called the \emph{support set}, e.g., in \cite[Section~1.7.1]{schneider2013convex})
is defined as follows \cite[Definition~3.1.3]{hiriart2012fundamentals}.
\begin{mydef}[{Exposed face}]
  \label{def:exposed}
  Consider a nonempty closed convex set $\Csetvar$.
  Given a vector $y \neq 0$, the \emph{exposed face} of $\Csetvar$ associated
  to $y$ is
\[
    F_\Csetvar(y) = \{\, x \in \Csetvar \mid \la x, y \ra = \supfun{y}{\Csetvar} \,\}.
\]
\end{mydef}

The exposed faces and normal cones are related by the following property \cite[Proposition~C.3.1.4]{hiriart2012fundamentals}.

\begin{myprop}
  \label{prop:normalsupfun}
  Consider a nonempty closed convex set $\Csetvar$.
  For any $x \in \Csetvar$ and nonzero vector $y$,
  $x \in F_\Csetvar(y)$ if and only if $y \in N_\Csetvar(x)$.
\end{myprop}

Given a set $\Sset$ and a matrix $A$,
let $A^{-\Tr}$ denote the preimage:
\begin{equation}
  \label{eq:preimage}
  A^{-\Tr}\Sset \eqdef \setbuild{x}{A^\Tr x \in \Sset}.
\end{equation}

\begin{myprop}[{\cite[Corollary~16.3.2]{rockafellar2015convex}}]
  \label{prop:podu}
  $ $\\For any convex set $\Csetvar$ and linear map $A$,
\[
    \polar{(A\Csetvar)} = A^{-\Tr} \polar{\Csetvar}.
\]
  where $\polar{\Csetvar}$ denotes the polar of the set $\Csetvar$.
\end{myprop}

\begin{myprop}[{\cite[Corollary~11.24(c)]{rockafellarvariational} or \cite[Corollary~16.3.1]{rockafellar2015convex}}]
  \label{prop:hA}
  Given a matrix $A \in \R^{n_1 \times n_2}$ and a nonempty closed convex
  set $\Sset \subseteq \R^{n_2}$, for all $y \in \R^{n_1}$,
  the following holds:
  \begin{equation}
    \label{eq:hA}
    \supfun{y}{A\Sset} = \supfun{A^\Tr y}{\Sset}.
  \end{equation}
\end{myprop}

\begin{myprop}[{\cite[Corollary~13.1.1]{rockafellar2015convex}}]
  \label{prop:supportincl}
  Consider two nonempty closed convex subsets $\Csetvar_1, \Csetvar_2 \subseteq \R^n$.
  The inclusion $\Csetvar_1 \subseteq \Csetvar_2$ is equivalent to the inequality
  $\supfun{x}{\Csetvar_1} \le \supfun{x}{\Csetvar_2}$ for all $x \in \R^n$.
\end{myprop}

When the support function is differentiable at a given point, $F_\Csetvar$ is a singleton and may be directly obtained using the following result:
\begin{myprop}[{\cite[Corollary~25.1.2]{rockafellar2015convex}}]
  \label{prop:exposed}
  $ $\\Given a nonempty closed convex set $\Csetvar$,
  if $\supfun{y}{\Csetvar}$ is differentiable at $y$ then $F_{\Csetvar}(y) = \{\grad \supfun{y}{\Csetvar}\}$.
\end{myprop}

In fact,
for nonempty compact convex sets,
the differentiability at $y$ is even
a necessary and sufficient conditions for
the uniqueness of $F_{\Csetvar}(y)$~\cite[Corollary~1.7.3]{schneider2013convex}.

\end{document}